\documentclass[12pt]{article}
\usepackage{amsmath}
\usepackage{amsthm}
\usepackage{amsfonts}
\usepackage{amssymb}

\newtheorem{thm}{Theorem}[section]
\newtheorem{lem}[thm]{Lemma}

\theoremstyle{definition}

\theoremstyle{remark}

\begin{document}
\title{On some properties of $\sigma(N)$
}
\author{Tomohiro Yamada}
\date{}
\maketitle

\begin{abstract}
We show asymptotic upper and lower bounds
for the greatest common divisor of $N$ and $\sigma(N)$.
We also show that there are infinitely many integers $N$
with fairly large g.c.d. of $N$ and $\sigma(N)$.
\end{abstract}

\section{Introduction}\label{intro}
  We denote by $\sigma(N)$ the sum of divisors of $N$.
  $N$ is said to be perfect if $\sigma(N)=2N$ and multiperfect
if $\sigma(N)=kN$ for some integer $k$.  It is not known whether or not
an odd perfect/multiperfect number exists.  There are known many results
which must be satisfied by such a number.  But these results are far from
solving whether or not an odd perfect/multiperfect number exists.

  Instead, we consider an analog of perfect/multiperfect numbers.  By definition,
we can easily see that $N$ is multiperfect if and only if $(N, \sigma(N))=N$.
  On the other hand, it is clear that $(N, \sigma(N))\leq N$ for any integer $N$.
  So we can suppose an integer $N$ is close to being perfect if $(N, \sigma(N))$ is large
relatively to $N$.

  Since $(N, \sigma(N))=1$ and $(N(N+1), \sigma(N(N+1))\geq N\geq (N(N+1))^{1/2}-1$
if $N$ is prime, we find that the behavior of $(N, \sigma(N))$ is very irregular.

  However, if we suppose that $N$ is out of some set of density zero,
we can obtain some nontrivial estimates.  In this direction,
Erd{\H o}s\cite{Erd} shows that $(N, \varphi(N))>1$ for almost all $N$.
His proof can be easily modified to prove that $(N, \sigma(N))>1$
for almost all $N$.  K\'atai and Subbarao\cite{KS} shows that for any given integer $l$,
the inequality $(N, \sigma(N+l))>1$ holds for almost all $N$.
  Before stating our result, we note that our arithmetic function does not
need to be the sum-of-divisors function.  Indeed, we only require that $\sigma$ is
multiplicative, $\sigma(p)=p+1$ for any prime $p$, and $\sigma(p^e)<p^{O(e)}$
uniformly for $p, e$ with $p$ prime.  For example, the sum of unitary divisors
function $\sigma^{*}(N)$ still works($d$ is called to be an unitary divisor of $N$
if $d$ divides $N$ and $(d, N/d)=1$.  Hence $\sigma^{*}(N)$ is multiplicative
and $\sigma^{*}(p^e)=p^e+1$).
\begin{thm}\label{th11}
  If $f(N)$ tends to infinity as $N$ does, then we have
\begin{equation}
(N, \sigma(N))\leq (\log \log N)^{f(N)}
\end{equation}
out of a set(depending on $f$) of density zero.
\end{thm}
\begin{thm}\label{th12}
  There exists some constant $c$ such that
\begin{equation}
(N, \sigma(N))\geq (\log \log \log N)^c
\end{equation}
holds for almost all integers $N$.
\end{thm}

  In these two theorems, $N$ runs over all integers.  We can see that
Theorem \ref{th11} remains valid if $N$ runs over shifted primes.
\begin{thm}\label{th13}
  If $f(N)$ tends to infinity as $N$ does, then we have
\begin{equation}
(p+a, \sigma(p+a))\leq (\log \log (p+a))^{f(p+a)}
\end{equation}
for any prime $p$ out of a set(depending on $f$ and $a$) of density zero.
\end{thm}

  Above three theorems are results valid for almost all members.
  To a contrary direction, it is interesting to find a function $f(N)$
which grows relatively fast such that $(N, \sigma(N))>f(N)$ for infinitely
many integers $N$.  As mentioned before, $f(N)=N^{1/2}-1$ works.
Using Theorem \ref{th13}, we can prove the following result.
\begin{thm}\label{th14}
  For any $\epsilon>0$, there exist infinitely many integers $N$ such that
\begin{equation}
(N, \sigma(N))\geq N^{2/3-\epsilon}.
\end{equation}
\end{thm}

\section{Notations and Preliminary Lemmas}\label{lemmas}
  Many of our notations and lemmas are due to
Katai and Wijsmuller \cite{KW}.  We denote by $P(n), p(n)$
the largest and smallest prime factor of $n$ respectively.
  For the variable $x$, let $x_1=\log x$, $x_2=\log x_1, \cdots$.
  We denote by $c$ some positive constant not necessarily
same at every occurence.  Let
\begin{equation}
s(x, k)=\sum_{p\leq x, p\equiv -1\pmod{k}}\frac{1}{p}.
\end{equation}
We cite some lemmas from \cite{KW}.  Lemmas \ref{lm21} and \ref{lm22}
are Lemma 1 in \cite{KW}.  Lemma \ref{lm23} is Lemma 2 in \cite{KW}.
\begin{lem}\label{lm21}
Let $Y_x$ and $E_x$ tend to infinity as $x$ does.  Write $n=n_1 n_2$
for each integer $n$ such that $P(n_1)\leq Y_x<p(n_2)$.  Let
\begin{equation*}
S_1(x)=\{n\mid n\leq x, n_1\geq Y_x^{E_x} \}.
\end{equation*}
  Then $\#S_1(x)/x$ tends to zero as $x$ goes to infinity.
\end{lem}
\begin{lem}\label{lm22}
  Let $Y_x$ and $E_x$ tend to infinity as $x$ does.  Write $n=n_1 n_2$
in the same way as in Lemma \ref{lm21}.  Let
\begin{equation*}
S_2(x)=\{n\mid n\leq x, n_2\text{ is not square free}\}.
\end{equation*}
  Then $\#S_2(x)/x$ tends to zero as $x$ tends to infinity.
\end{lem}
\begin{lem}\label{lm23}
  Uniformly in $k$ and $x\geq e^2$, we have
\begin{equation*}
s(x, k)\ll\frac{x_2}{\varphi(k)}.
\end{equation*}
\end{lem}

For the purpose to prove Theorem \ref{th13}, we need corresponding results
of these lemmas concerning the set of shifted primes.
\begin{lem}\label{lm24}
If $E_x$, $Y_x$, $E_x/\log\log Y_x$, $x/Y_x$ tend to infinity as $x$ does, we have
\begin{equation}
\sum_{y\geq Y_x^{E_x}, P(y)<Y_x}\pi(x, y, -a)=o(\pi(x)).
\end{equation}
\end{lem}
\begin{proof}
We divide the sum into two parts according to whether $y$ is large or small.

For $q<x^{1/2}$, by the Brun-Titchmarsh theorem, we have
\begin{equation}
\begin{split}
\sum_{Y_x^{E_x}\leq y<x^{1/2}, P(y)<Y_x}\pi(x, y, -a)
&\ll\frac{x}{x_1}\sum_{Y_x^{E_x}\leq y<x^{1/2}, P(y)<Y_x}\frac{\log x}{\varphi(y)\log(x/y)}\\
&\ll\frac{x}{x_1}\sum_{Y_x^{E_x}\leq y<x^{1/2}, P(y)<Y_x}\frac{1}{\varphi(y)}.
\end{split}
\end{equation}

It is well known that the number of integers $\leq X$ with largest prime factor $\leq Y$
is $\ll X\exp(-(\log X)/(\log Y))$.  It is also well known that $\varphi(n) \gg n/\log\log n$.
So partial summation gives
\begin{equation}\label{eq21}
\begin{split}
\sum_{Y_x^{E_x}\leq y<x^{1/2}, P(y)<E_x}\frac{1}{\varphi(y)}
&\ll\log\log(Y_x^{E_x})\exp(-\log(Y_x^{E_x})/\log Y_x)\\
&+\int_{Y_x^{E_x}}^{x^{1/2}}\frac{\log\log t}{t\exp(\frac{\log t}{\log Y_x})}dt.
\end{split}
\end{equation}
The first term is
\begin{equation}\label{eq22}
\leq\frac{\log E_x+\log\log Y_x}{e^{E_x}}\to 0
\end{equation}
as $x$ tends to infinity since $(\log\log Y_x)/E_x\to 0$ as $x\to\infty$.

For the same reason, we have
\begin{equation}
\log\log\log t/\log t\leq\log\log(E_x\log Y_x)/(E_x\log Y_x)=o(1/\log Y_x)
\end{equation}
and therefore we can majorize $\log\log t$ by $c\exp(\log t/(2\log Y_x))$.
Now the above integration is
\begin{equation}\label{eq23}
\begin{split}
\leq\int_{Y_x^{E_x}}^{x^{1/2}}\frac{1}{t\exp(\frac{\log t}{2\log Y_x})}dt
&=\int_{Y_x^{E_x}}^{x^{1/2}}t^{-(1+\frac{1}{2\log Y_x})}dt\\
&\ll (\log Y_x)Y_x^{-E_x/(2\log Y_x)}\\
&\ll\frac{\log Y_x}{\exp cE_x}\to 0
\end{split}
\end{equation}
as $x$ tends to infinity.  Substituting estimations (\ref{eq22}) and (\ref{eq23})
into (\ref{eq21}), we have
\begin{equation}\label{eq24}
\sum_{Y_x^{E_x}\leq y<x^{1/2}, P(y)<Y_x}\pi(x, y, -a)=o(x/x_1)=o(\pi(x)).
\end{equation}

For $y\geq x^{1/2}$, we obtain a trivial estimate
\begin{equation}\label{eq25}
\sum_{y\geq x^{1/2}, P(y)<Y_x}\pi(x, y, -a)\ll x\sum_{y\geq x^{1/2}, P(y)<Y_x}\frac{1}{y}.
\end{equation}
Using partial summation in a similar way to the case $y<x^{1/2}$, we have
\begin{equation}\label{eq26}
\begin{split}
\sum_{y\geq x^{1/2}, P(y)<Y_x}\frac{1}{y}&\leq\exp (-\log x/(2\log Y_x))+c\int_{x^{1/2}}^{\infty}\frac{dt}{t\exp(\log t/\log Y_x)}\\
&\ll\left(x^{1/(2\log Y_x)}+\int_{x^{1/2}}^{\infty}t^{-(1+\frac{1}{\log Y_x})}dt\right)\\
&\ll\frac{\log Y_x}{x^{1/(2\log Y_x)}}=o(\frac{x}{x_1^2})
\end{split}
\end{equation}
as $x$ tends to infinity since $\log Y_x=o(\log x)$ by assumption.  Substituting (\ref{eq26}) into (\ref{eq25}), we have
\begin{equation}\label{eq27}
\sum_{y\geq x^{1/2}, P(y)<Y_x}\pi(x, y, -a)=o(x/x_1)=o(\pi(x)).
\end{equation}
Combining (\ref{eq24}) and (\ref{eq27}) immediately gives the estimate in the lemma.
\end{proof}
\begin{lem}\label{lm25}
If $Y_x\to\infty$ as $x\to\infty$, then
\begin{equation}
\sum_{q\geq Y_x}\pi(x, q^2, -a)=o(\pi(x)).
\end{equation}
\end{lem}
\begin{proof}
As in the proof of Lemma \ref{lm24}, we divide the sum
into two parts.  For $q\leq x^{1/4}$, by the Brun-Titchmarsh theorem, we have
\begin{equation}\label{eq28}
\begin{split}
\sum_{Y_x\leq q<x^{1/4}}\pi(x, q^2, -a)
&\ll\frac{x}{x_1}\sum_{Y_x\leq q<x^{1/4}}\frac{x_1}{q(q-1)\log(x/q)}\\
&\ll\frac{x}{x_1}\sum_{Y_x\leq q<x^{1/4}}\frac{x_1}{q^2}\\
&\ll\frac{x}{Y_x x_1}=o(\frac{x}{x_1})=o(\pi(x))
\end{split}
\end{equation}
since $Y_x$ tends to infinity together with $x$.

For $q\geq x^{1/4}$, we obtain a trivial estimate
\begin{equation}\label{eq29}
\sum_{q\geq x^{1/4}}\pi(x, q^2, -a)\ll x\sum_{q\geq x^{1/4}}q^{-2}
=O(x^{3/4})=o(\pi(x)).
\end{equation}
Combining (\ref{eq28}) and (\ref{eq29}) immediately gives the stated inequality.
\end{proof}

\begin{lem}\label{lm26}
Let
\begin{equation}
t(x, q, a)=\sum_{p\leq x, p\equiv -1\pmod{q}}\frac{\pi(x, p, -a)}{\pi(x+a)}.
\end{equation}
Then we have
\begin{equation}
t(x, q, a)\ll\frac{x_2}{q}
\end{equation}
uniformly for $q\leq\log x$.
\end{lem}
\begin{proof}
As in previous lemmas, we divide the sum into two parts.  For $p\leq x/e$, by the Brun-Titchmarsh theorem, we have
\begin{equation}\label{eq2a}
\sum_{p\leq x/e, p\equiv -1\pmod{q}}\pi(x, p, -a)\leq\frac{cx}{x_1}\sum_{p\leq x/e, p\equiv -1\pmod{q}}\frac{x_1}{(p-1)\log(x/p)}.
\end{equation}
Let $\beta=\log p/\log x$.  Then we have $x_1/\log(x/p)=1/(1-\beta)$ and by partial summation,
\begin{equation}
\begin{split}
&\sum_{p\leq x/e, p\equiv -1\pmod{q}}\frac{x_1}{(p-1)\log(x/p)}\\
=&\sum_{p\leq x/e, p\equiv -1\pmod{q}}\frac{1}{(p-1)(1-\beta)}\\
\leq&\frac{\pi(x, q, -1)\log x}{x}+\int_2^{x/e}\frac{\pi(x, q, -1)}{t^2(1-\frac{\log t}{\log x})}dt\\
\ll&\frac{1}{q}\left(1+\int_2^{x/e}\frac{dt}{t\log t(1-\frac{\log t}{\log x})}\right).
\end{split}
\end{equation}
Setting $u=\log t$, the last integration can be computed and estimated as follows:
\begin{equation}
\begin{split}
&\int_2^{x/e}\frac{dt}{t\log t(1-(\log t)/(\log x))}\\
=&\int_{\log 2}^{\log(x/e)}\frac{du}{u(1-(u/\log x))}\\
=&\int_{\log 2}^{\log(x/e)}\frac{1}{u}+\frac{\log x}{1-(u/\log x)}du\\
=&\log\frac{\log(x/e)}{\log 2}+\log(1-\frac{\log 2}{\log x})-\log(1-\frac{(\log x)-1}{\log x})\\
\ll&\log\log x.
\end{split}
\end{equation}
Thus we can bound (\ref{eq2a}) by
\begin{equation}\label{eq2b}
\begin{split}
\sum_{p\leq x/e, p\equiv -1\pmod{q}}\pi(x, p, -a)\ll\frac{x x_2}{qx_1}.
\end{split}
\end{equation}

For $p\leq x/e$, we obtain a trivial estimate
\begin{equation}\label{eq2c}
\sum_{x/e\leq p\leq x, p\equiv -1\pmod{q}}\pi(x, p, -a)
\leq x\sum_{x/e\leq p\leq x, p\equiv -1\pmod{q}}\frac{1}{p}
=O(\frac{x}{qx_1})
\end{equation}
observing that $\log\log x-\log\log(x/e)=1/x_1$.
Combining (\ref{eq2b}) and (\ref{eq2c}), we obtain the stated inequality.
\end{proof}

\section{Proof of Theorem \ref{th11}}
Let $g(N)$ be an arbitrary function
tending to infinity.  We shall show that
the number of integers $n\leq x$ such that
$(n, \sigma(n))>(\log \log x)^{g(x)}$ is $o(x)$.
We can easily derive the theorem from this statement.
Let $g(x)=f(\exp \exp (\log \log x)^{1/2})/2$.
Then  $(\log \log x)^{g(x)}<(\log \log n)^{f(n)}$
if $\exp \exp n>(\exp \exp x)^{1/2}$.  Hence
the number of integers $n\leq x$ such that
$(n, \sigma(n))>(\log \log n)^{f(n)}$ is at most $o(x)+(\exp \exp x)^{1/2}=o(x)$.
Hence the theorem is proved.

By Lemmas \ref{lm21} and \ref{lm22},
we have $\#S_1(x)=o(x)$ and $\#S_2(x)=o(x)$.
So we may assume $n$ belongs to none of these sets.

We set $Y_x=E_x=x_4$.  Then $\sigma(p^e)<x_2$ if $p\leq Y_x$ and $e\leq E_x$.
Let $q\geq x_2$ be a prime dividing $(n, \sigma(n))$.
Then there exists an prime power divisor $p^e$ such that
$p^e\mid\mid n$ and $q\mid\sigma(p^e)$.
If $p\leq Y_x$, $p^e\leq Y_x^{E_x}$ since $n$ does not belong to $S_1$.
Hence $\sigma(p^e)<x_2\leq q$.  Thus $p$ must be greater than $Y_x$.
Now, since $n$ does not belong to $S_2$, $p^2$ does not divide $n$.
Hence $e=1$ and we have $p\equiv -1\pmod{q}$.  Therefore Lemma \ref{lm23}
gives that the number of $n\leq x$ such that $q\mid(n, \sigma(n))$
is bounded by
\begin{equation}
x\sum_{p\leq x, p\equiv -1\pmod{q}}\frac{1}{qp}=\frac{xs(x, q)}{q}
\ll x\frac{x_2}{q^2}.
\end{equation}
Hence we find that the number of $n\leq x$ such that $q$ divides
$(n, \sigma(n))$ for some $q\geq x_2$ is at most
\begin{equation}
cx\sum_{q\geq x_2}\frac{x_2}{q^2}\ll\frac{xx_2}{x_2\log x_2}=O\left(\frac{x}{x_3}\right).
\end{equation}

It follows that $(n, \sigma(n))$ divides $n_3$ with at most $o(x)$ exceptions,
where $n=n_3 n_4$ with $P(n_3)<x_2<p(n_4)$.  Our statement follows by observing
that the number of integers $n\leq x$ with $n_3\geq x_2^{g(x)}$ is at most $o(x)$
by Lemma \ref{lm21}.  This completes the proof of Theorem \ref{th11}.

\section{Proof of Theorem \ref{th12}}
Denote by $N(x, k, Q)$ the number of integers $n\leq x$ such that $n$ is divisible by at most $k$ primes $\leq Q$.  Brun's sieve gives $N(x, 0, Q)=O(x/\log Q)$ uniformly
for $Q, x$ with $Q<x$ and a simple inductive argument immediately gives
\begin{equation}\label{eq31}
N(x, k, Q)=O(x(c+\log\log Q)^k/\log Q)
\end{equation}
uniformly for $Q, x, k$ with $Q^{k+1}<x$.

Let $Q$ and $l$ be integers satisfying $x>Q^{l+1}$.

Let $q_1<q_2<\cdots<q_l<Q$ be distinct primes and
$M(y, q_1, \cdots, q_l)$ denote the number of integers $n<y$
such that $\sigma(n)$ is divisible by none of $q_1, \cdots, q_l$.
If $p$ is a prime $\equiv -1\pmod{q}$, then $p$
does not divide $n$ or $p^2$ divides $n$.  Hence it follows from
Brun's sieve and the Prime Number Theorem in arithmetic progressions
in the form of Theorem 9.6 in Karatsuba \cite{Krt} that
\begin{equation}\label{eq32}
M(y, q_1, \cdots, q_l)=O\left(\frac{ly}{(\log y)^{1/(Q-1)}}\right)
\end{equation}
uniformly for $Q, y$ with $Q<2\log y$.

Let $q_1<q_2<\cdots<q_l<Q$ be distinct primes.  Let $R_{xj}$ be functions of $x$
such that $R_{xj}<x^{1/2l}$.  Then the number of integers $n$
such that $n<x$, $q_1\cdots q_l\mid n$ and $q_1\cdots q_l\nmid\sigma(n)$ is at most

{\large
\begin{equation}
\begin{split}
&\sum_{e_1, \cdots, e_l\geq 1}M\left(\frac{x}{q_1^{e_1}\cdots q_l^{e_l}}, q_1, \cdots, q_l\right)\\
\leq &\sum_{e_1, \cdots, e_l\geq 1, \forall j, q_j^{e_j}\leq R_{xj}}M\left(\frac{x}{q_1^{e_1}\cdots q_l^{e_l}}, q_1, \cdots, q_l\right)\\
&+\sum_{j=1}^{l}\sum_{e_1, \cdots, e_l\geq 1, q_j^{e_j}\geq R_{xj}}\frac{x}{q_1^{e_1}\cdots q_l^{e_l}}\\
&=\Sigma_1+\Sigma_2, \text{ say}.
\end{split}
\end{equation}
}

By (\ref{eq32}), we see that $\Sigma_1$ is
\begin{equation}
O\left(\frac{lx}{(q_1-1)\cdots (q_l-1) x_1^{1/(Q-1)}}\right)
\end{equation}
since $Q<2\log(x^{1/2})<\log(x/(R_{x1}\cdots R_{xl}))<\log(x/(q_1^{e_1}\cdots q_l^{e_l}))$.
A trivial argument gives
\begin{equation}
\Sigma_2\leq\sum_{j=1}^{l}\frac{xq_j}{R_{xj}(q_1-1)\cdots(q_l-1)}\ll lx^{1-1/2l}.
\end{equation}
Combining these estimates, we obtain
\begin{equation}
\begin{split}
&\sum_{e_1, \cdots, e_l\geq 1}M\left(\frac{x}{q_1^{e_1}\cdots q_l^{e_l}}, q_1, \cdots, q_l\right)\\
\ll&\left(\frac{lx}{(q_1-1)\cdots (q_l-1) x_1^{1/(Q-1)}}+lx^{1-1/(2l)}\right).
\end{split}
\end{equation}
Hence the number of integers $n<x$ such that there exist distinct primes
$q_1<q_2<\cdots<q_l<Q$ satisfying $q_1\cdots q_l\mid n$ and $q_1\cdots q_l\nmid\sigma(n)$
is bounded by
\begin{equation}\label{eq33}
\begin{split}
&\sum_{q_1<q_2<\cdots<q_l<Q}c\left(\frac{lx}{(q_1-1)\cdots (q_l-1) x_1^{1/(Q-1)}}+lx^{1-1/2l}\right)\\
\ll&\frac{lx(\log\log Q+c)^l}{x_1^{1/(Q-1)}}+\frac{lQ^lx}{x^{1/2l}}.
\end{split}
\end{equation}

We observe that $Q=x_2/x_3$ and $l=cx_4/x_5$ with $c>0$ sufficiently small satisfy
our conditions and the right-hand side of (\ref{eq33}) is $o(x)$.
By (\ref{eq31}), the number of integers $n<x$ divisible by at most $l-1$
distinct primes smaller than $Q$ is $o(x)$.
Now the remaining integers $n$ have the property that $(n, \sigma(n))$ has
at least $l$ distinct prime factors.  Hence $(n, \sigma(n))>l^{cl}>x_3^c$.
This completes the proof.

\section{Proof of Theorem \ref{th13} and \ref{th14}}
Set $Y_x=E_x=x_4$ as in the proof of Theorem \ref{th11}.
We easily see that this choice satisfies the condition of Lemma \ref{lm24}.
Proceeding in the same way as our proof of Theorem \ref{th11},
we immediately obtain Theorem \ref{th13} using Lemmas \ref{lm24}-\ref{lm26}
instead of Lemmas \ref{lm21}-\ref{lm23}.

Let $N=p(p+1)m$ where $m$ is the largest divisor of $\sigma(p+1)$
relatively prime to $p+1$.  The proof of Theorem \ref{th13} shows that
$P((p+1, \sigma(p+1)))\leq n_1$ and $n_1<(\log\log(p+1))^{f(p+1)}$
for almost all prime $p$, where $p+1$ is decomposed into $n_1 n_2$ such that $P(n_1)\leq\log\log(p+1)<p(n_2)$.  Hence if $\epsilon$ is an arbitrary positive real number,
$(p+1, \sigma(p+1))=o(x^{\epsilon})$ holds for almost all prime $p$.

Hence there exists infinitely many prime $p$ such that $m>p^{1-\epsilon}$.
If we choose such $p$, then $(p+1)m$ divides $\sigma(N)$.
Hence $(p+1)m$ divides $(N, \sigma(N))$ and clearly $(p+1)m>N^{2/3-\epsilon}$.
This completes the proof.

{}
\vskip 12pt

{\small Tomohiro Yamada}\\
{\small Department of Mathematics\\Faculty of Science\\Kyoto University\\Kyoto, 606-8502\\Japan}\\
{\small e-mail: \protect\normalfont\ttfamily{tyamada@math.kyoto-u.ac.jp}}
\end{document}